\documentclass{article}

\usepackage{arxiv}

\usepackage[utf8]{inputenc} 
\usepackage[T1]{fontenc}   
\usepackage{hyperref}    
\usepackage{url}           
\usepackage{booktabs}      
\usepackage{amsfonts}       
\usepackage{nicefrac}    
\usepackage{microtype}      
\usepackage{lipsum}
\usepackage{changepage}
\usepackage{graphicx}
\usepackage{amsmath}
\usepackage{amsfonts}
\usepackage{amssymb}
\usepackage{makeidx}
\usepackage{algorithm,algorithmic}
\usepackage[switch]{lineno}
\usepackage{xcolor}
\usepackage{bm}
\usepackage[justification=centering]{caption}
\usepackage{amsfonts}
\usepackage[bbgreekl]{mathbbol}
\DeclareSymbolFontAlphabet{\mathbbm}{bbold}
\DeclareSymbolFontAlphabet{\mathbb}{AMSb}%
\usepackage{amsthm}
\newtheorem{problem}{Problem}

\title{Model Predictive Control using MATLAB}

\author{
  Midhun T. Augustine\\
  Department of Electrical Engineering\\
  Indian Institute of Technology Delhi, India \\
  midhunta30@gmail.com \\
 \vspace{.01cm}\\
  Date of initial version: 15 - 12 - 2021\\ Date of current version: 01 - 09 - 2023\\ \vspace{.01cm}}

\begin{document}
\maketitle

\begin{abstract}
This tutorial consists of a brief introduction to the modern control approach called model predictive control (MPC) and its numerical implementation using MATLAB. We discuss the basic concepts and numerical implementation of the two major classes of MPC: Linear MPC (LMPC) and Nonlinear MPC (NMPC). This includes the various aspects of MPC such as formulating the optimization problem,  constraints handling, feasibility, stability, and optimality.
\end{abstract}


\section{Introduction}
MPC is a feedback control approach that uses model-based optimization for computing the control input.
In MPC, a model of the system along with the current state (measured or estimated) 
is used to predict the future behavior (states) of the system, for a control input sequence over a short period. The predicted behavior is characterized by a cost function which is a function of the predicted state and control sequence.
Then an optimization algorithm is used to find the control sequence which optimizes the predicted behavior or cost function.
The first element of the control sequence is applied to the system which gives the next state, and the algorithm is repeated at the next time instant, which results in a receding horizon scheme. Model predictive control (MPC) is also known as receding horizon control (RHC).
The name MPC originated from the model-based predictions used for optimization, whereas the name RHC comes from the receding horizon nature of the control scheme. MPC which originated from the optimal control approach has the following advantages over the former:
\begin{enumerate}
    \item It gives closed-loop control schemes whereas optimal control mostly results in open-loop control schemes.
    \item MPC can handle complex systems such as nonlinear, higher-order, multi-variable, etc. 
   \item MPC can incorporate constraints easily.
\end{enumerate}
\par \textit{Notations:}
$\mathbb{N},\mathbb{Z}$ and $\mathbb{R}$ denote the set of natural numbers, integers, and real numbers respectively.
$\mathbb{R}^{n}$ stands for $n$ - dimensional Euclidean space and $\mathbb{R}^{m \times n}$ refers to the space of $m \times n$ real matrices.  
Matrices and vectors are represented by boldface letters ($\textbf{A},\textbf{a}$),  scalars by normal font ($A,a$), and sets by blackboard bold font ($\mathbb{A},\mathbb{B},\mathbb{C}$). The notation $\textbf{P}>0$ ($\textbf{P}\geq 0$) indicates that $\textbf{P}$ is a real symmetric positive definite (semidefinite) matrix. Finally, $\textbf{I},\textbf{0}$   represents the identity matrix and zero matrix of appropriate order.
\begin{figure} 
 		\begin{center}
 		\includegraphics [scale=.375] {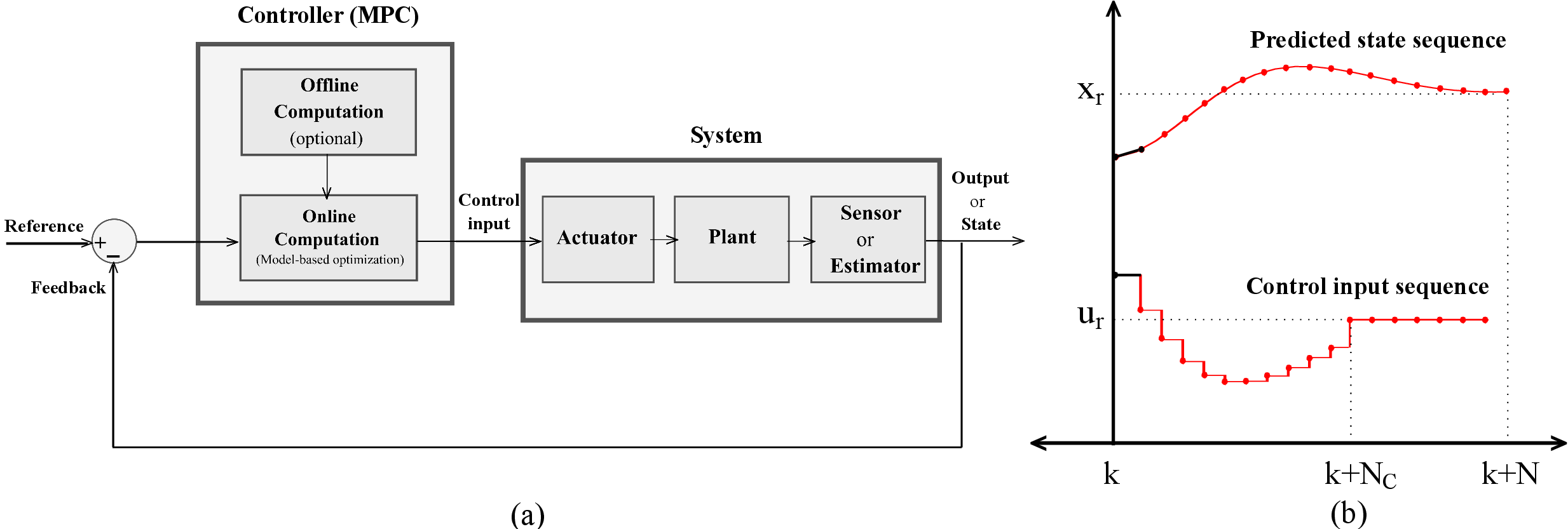}
 		\caption{{\footnotesize (a) MPC  General block diagram \hspace{.2cm} (b) Basic MPC  strategy.}}
 	\end{center}
 \end{figure}
\par MPC is associated with a number of terminologies which are defined below
\begin{enumerate}
        \item \textbf{Sampling time} (T): It is the time difference between two consecutive state measurements or control updates. In general $T\in \mathbb{R}^{+},$ and for discrete-time systems $T>0,$ whereas for continuous-time systems $T=0.$

    \item \textbf{Time horizon} ($N_{T}$): It is the number of time instants the control input is applied to the system. In general, $N_T \in \mathbb{N}$ and if $N_T$ is infinity, the problem is called an infinite horizon problem, otherwise finite horizon problem.
    \item \textbf{Prediction horizon} ($N$): It is the length of the prediction window over which the states are predicted and optimized. In general $N\in \mathbb{N}$ and usually $2\leq N\leq N_{T}.$
    \item \textbf{Control horizon} ($N_C$): It is the length of the control window in which the control input is optimized, and normally $N_C \leq N$. In this tutorial, we mainly focus on the case for which the control horizon is the same as the prediction horizon, i.e., $N_C=N.$
    If $N_C<N$ we optimize the control sequences over $N_C$ and the remaining control sequences (of length $N-N_C$) are normally chosen as zero. 
\end{enumerate}
The general block diagram of a dynamical system with MPC is given in Fig. 1(a) and
the basic strategy of MPC is given in Fig. 1(b). In MPC during the
current time instant $k,$ we consider the optimization over the next ${N}$ instants where
${N}$ is the prediction horizon, i.e., the optimization window is from $k$ to $k+{N}.$ This indicates that the optimization window moves with time and this feature is called moving horizon or receding horizon. 
In MPC, during every time instant, we compute the sequence of control inputs over the control horizon, which optimizes the future performance of the system over the prediction horizon. Then the first element of the optimal control sequence is applied to the system, which results in a receding horizon scheme. The first element of the control sequence and the next state under the MPC scheme are represented in black color in Fig. 1(b). By repeating this at each time instant, we obtain the control inputs and states with MPC over the time horizon.

\par Based on the nature of the system model used in the optimization, MPC can be grouped into the following two classes:
\begin{enumerate}
    \item \textbf{Linear MPC}: For which the system model and the constraints are linear. 
    The cost function can be linear or quadratic which results in linear programming or quadratic programming problems which are convex optimization problems.
    \item \textbf{Nonlinear MPC}: For which the system model is nonlinear and the constraints can be either linear or nonlinear. The cost function is usually chosen as a linear or quadratic function of states and control inputs which results in a nonlinear programming problem that can be non-convex.
\end{enumerate}
Another classification is based on the implementation of the MPC algorithm which results in the following categories:
\begin{enumerate}
    \item \textbf{Implicit MPC}: This is also known as the traditional MPC in which the control input at each time instant is computed by solving an optimization problem online. In this tutorial, we will be focusing on implicit MPC which is the most general MPC scheme.
    \item \textbf{Explicit MPC}: In this, the online computation is reduced by transferring the optimization problem offline. In explicit MPC the state constraint set is divided into a finite number of regions and the optimization problem is solved offline for each of the regions which gives the control input as a function of the state. This simplifies the online computation to just a function evaluation.
\end{enumerate}
\par When it comes to optimization the approaches can be classified into two categories:
\begin{enumerate}
    \item \textbf{Iterative approach}: In which the elements of the decision vector are optimized together. Here the optimal decision vector is computed iteratively by starting with an initial guess which is then improved in each iteration. Most of the linear programming and nonlinear programming algorithms are based on the iterative approach. 
    \item \textbf{Recursive approach}: In which the elements of the decision vector are optimized recursively, i.e., one at a time. The popular optimization algorithm which uses the recursive approach is dynamic programming. Even though both the iterative approach and recursive approach are used in MPC, in this tutorial we focus on the former.
\end{enumerate}
 \newpage
 \section{MPC of Linear Systems}
 In this section, we discuss the basic concept of Linear MPC (LMPC) and its numerical implementation.
\subsection{LMPC: Problem Formulation}
Consider the discrete-time linear time-invariant (LTI) system:
\begin{equation}
{ \textbf{x}_{{k+1}}}=\textbf{A}\textbf{x}_{k}+\textbf{B}\textbf{u}_{k}
\end{equation}
where $k\in \mathbb{T}=\{0,1,...,N_{T}-1\}$ is the discrete time instant,
$\textbf{x}_{k}\in \mathbb{X} \subseteq \mathbb{R}^{n}$ is the state vector,  $\textbf{u}_{k}\in \mathbb{U} \subseteq \mathbb{R}^{m}$ is the control input vector, $\textbf{A} \in \mathbb{R}^{n \times n}$ is the system matrix and
$\textbf{B} \in \mathbb{R}^{n \times m}$ is the input matrix. The sets $\mathbb{X}$ and $\mathbb{U}$ are the constraint sets for the states and control inputs which are usually represented by linear inequalities:
\begin{equation}
\begin{aligned}
  & \mathbb{X}=\{\textbf{x} \in \mathbb{R}^{n}:  \textbf{F}_{\textbf{x}}\textbf{x} \leq \textbf{g}_{\textbf{x}} \}\\
  &\mathbb{U}=\{\textbf{u} \in \mathbb{R}^{m}:  \textbf{F}_{\textbf{u}}\textbf{u} \leq \textbf{g}_{\textbf{u}} \}.
  \end{aligned}
\end{equation}

The cost function is chosen as  a quadratic sum of the states  and control inputs:
\begin{equation}
J=\textbf{x}_{N_T}^{T}\textbf{Q}_{N_T}\textbf{x}_{N_T}+ \sum_{k=0}^{{N_T}-1}\textbf{x}_{k}^{T}\textbf{Q}\textbf{x}_{k}+\textbf{u}_{k}^{T}\textbf{R}\textbf{u}_{k} 
\end{equation}
where $\textbf{Q}_{N_T}\in \mathbb{R}^{n\times n},\textbf{Q} \in \mathbb{R}^{n \times n},$ $\textbf{R}  \in \mathbb{R}^{m \times m}$ are the weighting matrices used for relatively weighting the states and control inputs and to be chosen such that $\textbf{Q}_{N_T} \geq 0,$ $\textbf{Q} >0,$ $\textbf{R}>0$. The state and control sequence  is defined as $\textbf{X}=\big(\textbf{x}_{0}, \textbf{x}_{1},...,\textbf{x}_{{N_T}} \big),$
$\textbf{U}= 
\big(\textbf{u}_{0},\textbf{u}_{1},..., \textbf{u}_{{N_T}-1}
\big)$
which contains the state and control input over the time horizon.
 Now, the optimal control problem for the LTI system is defined as follows which is also known as the constrained linear quadratic regulator (CLQR) problem:
\begin{problem}
For the linear system (1) with the initial state $\textbf{x}_{0}$,
compute the control sequence $\textbf{U}$  by solving the optimization problem
\begin{equation}
    \begin{aligned}
    \underset{\textbf{U}}{\mbox{inf}} ~&~  J\\ \text{subject to}
    ~&~ \textbf{U} \in \mathbb{U}^{N_T}, \hspace{.2cm} \textbf{X}\in \mathbb{X}^{{N_T}+1}\\
    ~&~ \textbf{x}_{k+1}=\textbf{A}\textbf{x}_{k}+\textbf{B}\textbf{u}_{k}, \hspace{.3cm}k\in \mathbb{T}.
    \end{aligned}
\end{equation}
\end{problem}
As $N_T \rightarrow \infty$ the problem is called infinite-horizon constrained LQR. One can solve constrained LQR with a large time horizon ($N_T \rightarrow \infty$) using the MPC approach which usually results in suboptimal solutions with lesser computation. MPC uses a prediction horizon ${N} \leq N_T$ (in practice ${N} << N_T$) and during every time instant the control sequence for the next ${N}$ instants is computed for minimizing the cost over the next ${N}$ instants.
The cost function for the MPC with a prediction horizon $N$ at time instant $k$ is defined as
\begin{equation}
    J_k=\textbf{x}_{k+{N}|k}^{T}\textbf{Q}_{{N}}\textbf{x}_{k+{N}|k}+\sum_{i=k}^{k+{N}-1}\textbf{x}_{i|k}^{T}\textbf{Q}\textbf{x}_{i|k}+\textbf{u}_{i|k}^{T}\textbf{R}\textbf{u}_{i|k} 
\end{equation}
in which $\textbf{x}_{i|k},\textbf{u}_{i|k}$ denotes the state and control input at time instant $i$ predicted or computed at time instant $k.$ Note that here $k$ denotes the time instants within the time horizon and $i$ denotes the time instants within the prediction horizon.
Similarly, the state and control sequence for the MPC at time instant $k$ is defined as $\textbf{X}_{k}=\big(\textbf{x}_{k|k}, \textbf{x}_{k+1|k},...,\textbf{x}_{k+{N}|k} \big),$
$\textbf{U}_{k} = \big(
\textbf{u}_{k|k}, \textbf{u}_{k+1|k},...,\textbf{u}_{k+{N}-1|k} \big).$
Then the MPC problem for linear systems is defined as follows:
\begin{problem}
For the linear system (1) with the current state $\textbf{x}_{k|k}=\textbf{x}_{k}$ given,
compute the control sequence $\textbf{U}_{k},$ by solving the optimization problem
\begin{equation}
    \begin{aligned}
    \underset{\textbf{U}_{k}}{\mbox{inf}} ~&~  J_{k}\\ 
    \text{subject to}
    ~&~ \textbf{U}_{k} \in \mathbb{U}^{N},\hspace{.2cm} \textbf{X}_{k}\in \mathbb{X}^{N+1}, \hspace{.4cm}k\in \mathbb{T}\\
    ~&~ \textbf{x}_{i+1|k}=\textbf{A}\textbf{x}_{i|k}+\textbf{B}\textbf{u}_{i|k}, \hspace{.3cm}k\in \mathbb{T}, i=k,...,k+N-1.
    \end{aligned}
\end{equation}
\end{problem}.
\newpage
\subsection{LMPC: Algorithm}
Here we represent the MPC  optimization problem as a quadratic programming problem. 
From the solution of the state equation for LTI systems  we obtain
 \begin{equation}
 \left[\begin{matrix}
\textbf{x}_{k|k}\\\textbf{x}_{k+1|k}\\\vdots\\ \textbf{x}_{k+{N}|k}
\end{matrix}\right]=    \left[\begin{matrix}
\textbf{I}\\\textbf{A}\\ \vdots \\ \textbf{A}^{{N}} 
\end{matrix}\right] \textbf{x}_{k}+ \left[\begin{matrix}
\textbf{0} & \textbf{0}  & \dots &\textbf{0}  \\\textbf{B}& \textbf{0}  & \dots &\textbf{0}\\ \vdots & \vdots& & \vdots\\ \textbf{A}^{{N}-1}\textbf{B} & \textbf{A}^{{N}-2}\textbf{B}  & \dots &\textbf{B}   
\end{matrix}\right]\left[\begin{matrix}
\textbf{u}_{k|k}\\\textbf{u}_{k+1|k}\\\vdots\\ \textbf{u}_{k+{N}-1|k}
\end{matrix}\right].
\end{equation}
By defining the following matrices  
 \begin{equation}
\textbf{X}_{k} = \left[\begin{matrix}
\textbf{x}_{k|k}\\\textbf{x}_{k+1|k}\\\vdots\\ \textbf{x}_{k+{N}|k}
\end{matrix}\right], \hspace{.2cm}\textbf{U}_{k} = \left[\begin{matrix}
\textbf{u}_{k|k}\\\textbf{u}_{k+1|k}\\\vdots\\ \textbf{u}_{k+{N}-1|k}
\end{matrix}\right]\hspace{.2cm} \textbf{A}_{\textbf{X}}=    \left[\begin{matrix}
\textbf{I}\\\textbf{A}\\ \vdots \\ \textbf{A}^{{N}} 
\end{matrix}\right], \textbf{B}_{\textbf{U}}=\left[\begin{matrix}
\textbf{0} & \textbf{0}  & \dots &\textbf{0}  \\\textbf{B}& \textbf{0}  & \dots &\textbf{0}\\ \vdots & \vdots& & \vdots\\ \textbf{A}^{{N}-1}\textbf{B} & \textbf{A}^{{N}-2}\textbf{B}  & \dots &\textbf{B}   
\end{matrix}\right] \hspace{.2cm}
\end{equation}
  the equation (7) is rewritten as
\begin{equation}
    \textbf{X}_{k}=\textbf{A}_{\textbf{X}}\textbf{x}_{k}+\textbf{B}_{\textbf{U}}\textbf{U}_{k}.
\end{equation}
This indicates that, the predicted state $\textbf{X}_{k}$ can be represented as a function of the current state $\textbf{x}_{k}$ and   input sequence $\textbf{U}_{k}.$
Similarly, by defining
\begin{equation}
\textbf{Q}_{\textbf{X}} = \left[\begin{matrix}
\textbf{Q} & \dots & \textbf{0} & \textbf{0}\\ \vdots &  &\vdots & \vdots \\ \textbf{0} & \dots & \textbf{Q} & \textbf{0} \\ \textbf{0}& \dots  &\textbf{0} & \textbf{Q}_{N}
\end{matrix}\right], \hspace{.2cm}\textbf{R}_{\textbf{U}} = \left[\begin{matrix}
\textbf{R} &\textbf{0} &\dots  &\textbf{0} \\\textbf{0} & \textbf{R}  & \dots& \textbf{0} \\\vdots & \vdots&  &\vdots \\\textbf{0}& \textbf{0} & \dots & \textbf{R}
\end{matrix}\right]
\end{equation}
 the cost function (5) can be represented in terms of $ \textbf{X}_{k}$ and  $\textbf{U}_{k}$ as
\begin{equation}
    J_{k}=  \textbf{X}_{k}^{T}\textbf{Q}_{\textbf{X}} \textbf{X}_{k}+\textbf{U}_{k}^{T}\textbf{R}_{\textbf{U}} \textbf{U}_{k}.
\end{equation}
Finally, by defining 
\begin{equation}
\textbf{F}_{\textbf{X}} = \left[\begin{matrix}
\textbf{F}_{\textbf{x}} &\textbf{0} & \dots & \textbf{0} \\\textbf{0} & \textbf{F}_{\textbf{x}} & \dots & \textbf{0} \\ \vdots &\vdots  &  & \vdots \\\textbf{0} & \textbf{0} & \dots & \textbf{F}_{\textbf{x}}
\end{matrix}\right], \hspace{.2cm}
\textbf{g}_{\textbf{X}} = \left[\begin{matrix}
\textbf{g}_{\textbf{x}}\\\textbf{g}_{\textbf{x}} \\\vdots\\\textbf{g}_{\textbf{x}}
\end{matrix}\right],
\hspace{.2cm}
\textbf{F}_{\textbf{U}} = \left[\begin{matrix}
\textbf{F}_{\textbf{u}} & \textbf{0} & \dots & \textbf{0}\\ \textbf{0}& \textbf{F}_{\textbf{u}} & \dots & \textbf{0} \\\vdots &\vdots  & & \vdots \\ \textbf{0}& \textbf{0}& \dots & \textbf{F}_{\textbf{u}}
\end{matrix}\right], \hspace{.2cm}
\textbf{g}_{\textbf{U}} = \left[\begin{matrix}
\textbf{g}_{\textbf{u}}\\\textbf{g}_{\textbf{u}} \\\vdots\\\textbf{g}_{\textbf{u}}
\end{matrix}\right]
\end{equation}
the state and control constraints in (2) can be represented in terms of $ \textbf{X}_{k}$ and $\textbf{U}_{k}$  as
\begin{equation}
\begin{aligned}
  & \textbf{F}_{\textbf{X}}\textbf{X}_{k} \leq \textbf{g}_{\textbf{X}}\\
   & \textbf{F}_{\textbf{U}}\textbf{U}_{k} \leq \textbf{g}_{\textbf{U}}. 
\end{aligned}    
\end{equation}
Now by combining $\textbf{X}_{k}$ and $\textbf{U}_{k},$ we can represent the cost function with a single decision vector. For that we define 
\begin{equation}
\textbf{z} = \left[\begin{matrix}
\textbf{X}_{\textbf{k}}\\ \textbf{U}_{\textbf{k}}
\end{matrix}\right], \hspace{.2cm}\textbf{H}= \left[\begin{matrix}
\textbf{Q}_{\textbf{X}} & \textbf{0} \\\textbf{0} & \textbf{R}_{\textbf{U}}
\end{matrix}\right],
\hspace{.2cm}\textbf{F} = \left[\begin{matrix}
\textbf{F}_{\textbf{X}} &\textbf{0} \\\textbf{0} & \textbf{F}_{\textbf{U}}
\end{matrix}\right], \hspace{.2cm}\textbf{g} = \left[\begin{matrix}
\textbf{g}_{\textbf{X}} \\\textbf{g}_{\textbf{U}}
\end{matrix}\right],\hspace{.2cm}\textbf{F}_{eq} = \left[\begin{matrix}
\textbf{I}& -\textbf{B}_{\textbf{U}}
\end{matrix}\right]
, \hspace{.2cm}\textbf{g}_{eq} =\textbf{A}_{\textbf{X}}\textbf{x}_{k}
\end{equation}
using this we can rewrite the cost function (11) and constraints (9),(13) and represent the optimization problem (6) as a quadratic programming problem as below
\begin{equation}
\begin{aligned}
 \underset{\textbf{z}}{\text{inf}} ~&~     \textbf{z}^{T}\textbf{H}\textbf{z}\\ 
 \text{subject to}
~&~ \textbf{F}\textbf{z} \leq \textbf{g}\\
  ~&~ \textbf{F}_{eq}\textbf{z} = \textbf{g}_{eq}
     \end{aligned}
\end{equation}
which can be solved using standard numerical optimization algorithms such as the steepest-descent method, Newton method, etc.
For faster convergence of the numerical optimization method, the optimal solution for the current instant can be used as the initial condition for the next instant.
Note that here $\textbf{g}_{eq}$ is a function of the state vector $\textbf{x}_{k}.$ Therefore the current state information is required for solving the optimization problem.  
In MPC this optimization problem is solved during each time instant $k$ and the first element of $\textbf{U}_{k}^{*}$ is applied to the system, i.e., the control input with MPC  is 
\begin{equation}
    \textbf{u}_{k}=[\textbf{U}_{k}^{*}]_{1}=\textbf{u}_{k|k}^{*}.
\end{equation}
Note that this algorithm is based on the assumption that, an optimal control sequence exists at each time instant. The existence of an optimal control sequence depends on the system model and constraints, and this will be discussed in the feasibility analysis section.
The algorithm for linear MPC is given below:
\begin{algorithm}[H]
 \small
	\begin{algorithmic}[1] 
	
	\STATE Require $\textbf{A},\textbf{B},N_{T},N,n,m,\textbf{Q},\textbf{R},\textbf{Q}_{N_T},\textbf{F}_{\textbf{x}},\textbf{g}_{\textbf{x}},\textbf{F}_{\textbf{u}},\textbf{g}_{\textbf{u}}$
		\STATE Initialize $\textbf{x}_{0},\textbf{z}_{0}$
		\STATE Construct $\textbf{A}_{\textbf{X}},\textbf{B}_{\textbf{U}},\textbf{Q}_{X},\textbf{R}_{\textbf{U}},\textbf{H},\textbf{F},\textbf{g}$
		\FOR  {$k= 0~to~ N_{T}-1 $}
		\STATE $\textbf{x}_{k}=[\textbf{X}]_{k+1}$ (obtain $\textbf{x}_{k}$ from measurement or estimation)
		\STATE Compute $\textbf{F}_{eq},\textbf{g}_{eq}$
		\STATE  Compute $\textbf{z}^{*}=\left[\begin{matrix}\textbf{X}_{k}^{*}\\\textbf{U}_{k}^{*}\end{matrix}\right]$ by solving the optimization problem (15)
		\STATE Apply  $\textbf{u}_{k}=[\textbf{U}_{k}^{*}]_{1}$ to the system
		\STATE Update $\textbf{z}_{0}=\textbf{z}^{*}$
		\ENDFOR
	\end{algorithmic}
	\caption{: LMPC}
\end{algorithm}
The optimization problem can be solved using the MATLAB function $\textbf{fmincon}$ for solving  constrained optimization problems which are of the form
    \begin{equation}
        \textbf{z}^{*}=\mbox{fmincon}(f,\textbf{z}_{0},\textbf{F},\textbf{g},\textbf{F}_{eq},\textbf{g}_{eq},\textbf{lb},\textbf{ub})
    \end{equation}
    in which $\textbf{lb},\textbf{ub}$ are the vectors containing the lower bound and upper bound of each element in the decision vector $\textbf{z}.$
\subsection{Reducing online computation}
Here we discuss some methods for reducing online computation in which the basic idea is to reduce the number of optimization variables and constraints. The first method uses the idea of eliminating the states from the decision vector $\textbf{z}$. This method is useful when we have only control constraints, i.e., the state is unconstrained $\textbf{x}_{k} \in \mathbb{R}^{n}$ or the state constraints can be transferred to control constraints. We have from (11) the cost $J_k$ is a function of the state sequence $\textbf{X}_{k}$ and control sequence $\textbf{U}_{k}.$ Now, by substituting (9) in (11), we obtain
\begin{equation}
\begin{aligned}
       J_{k}&= \big[\textbf{A}_{\textbf{X}}\textbf{x}_{k}+\textbf{B}_{\textbf{U}}\textbf{U}_{k}\big]^{T}\textbf{Q}_{\textbf{X}} \big[\textbf{A}_{\textbf{X}}\textbf{x}_{k}+\textbf{B}_{\textbf{U}}\textbf{U}_{k}\big]+\textbf{U}_{k}^{T}\textbf{R}_{\textbf{U}} \textbf{U}_{k}\\
       &=\textbf{U}_{k}^{T}\big[\textbf{B}_{\textbf{U}}^{T}\textbf{Q}_{\textbf{X}}\textbf{B}_{\textbf{U}}+\textbf{R}_{\textbf{U}}\big]\textbf{U}_{k}+2 \textbf{x}_{k}^{T}\big[\textbf{A}_{\textbf{X}}^{T}\textbf{Q}_{\textbf{X}}\textbf{B}_{\textbf{U}} \big]\textbf{U}_{k}+ \textbf{x}_{k}^{T}\big[\textbf{A}_{\textbf{X}}^{T}\textbf{Q}_{\textbf{X}}\textbf{A}_{\textbf{X}} \big]\textbf{x}_{k}\\
       &=\textbf{U}_{k}^{T}\textbf{H}\textbf{U}_{k}+ \textbf{q}_{k}^{T}\textbf{U}_{k}+{r}_{k}
\end{aligned}
\end{equation}
where $\textbf{H}=\textbf{B}_{\textbf{U}}^{T}\textbf{Q}_{\textbf{X}}\textbf{B}_{\textbf{U}}+\textbf{R}_{\textbf{U}}, \textbf{q}_{k}^{T}=2\textbf{x}_{k}^{T}\textbf{A}_{\textbf{X}}^{T}\textbf{Q}_{\textbf{X}}\textbf{B}_{\textbf{U}}$ and ${r}_{k}=\textbf{x}_{k}^{T}\textbf{A}_{\textbf{X}}^{T}\textbf{Q}_{\textbf{X}}\textbf{A}_{\textbf{X}}\textbf{x}_{k}.$ Therefore we can represent the cost $J_k$ as a function of the current state $\textbf{x}_{k}$ and control sequence $\textbf{U}_{k},$ in which $\textbf{U}_{k}$ is the decision vector. Similarly, the constraint inequalities (13) can be rewritten as 
\begin{equation}
    \begin{aligned}
     & \textbf{F}_{\textbf{X}}\big[\textbf{A}_{\textbf{X}}\textbf{x}_{k}+\textbf{B}_{\textbf{U}}\textbf{U}_{k}\big] \leq \textbf{g}_{\textbf{X}} \implies  \textbf{F}_{\textbf{X}}\textbf{B}_{\textbf{U}}\textbf{U}_{k} \leq \textbf{g}_{\textbf{X}}-\textbf{F}_{\textbf{X}}\textbf{A}_{\textbf{X}}\textbf{x}_{k}\\
     &\textbf{F}_{\textbf{U}}\textbf{U}_{k} \leq \textbf{g}_{\textbf{U}}.
    \end{aligned}
\end{equation}
Now, by defining $\textbf{z}=\textbf{U}_{k},$ $\textbf{F}=\left[\begin{matrix}\textbf{F}_{\textbf{X}}\textbf{B}_{\textbf{U}}\\ \textbf{F}_{\textbf{U}}\end{matrix}\right],$ $\textbf{g}=\left[\begin{matrix}\textbf{g}_{\textbf{X}}-\textbf{F}_{\textbf{X}}\textbf{A}_{\textbf{X}}\textbf{x}_{k}\\ \textbf{g}_{\textbf{U}}\end{matrix}\right]$ we can represent the optimization problem (15) as a quadratic programming problem as below
\begin{equation}
\begin{aligned}
 \underset{\textbf{z}}{\text{inf}}    ~&~ \textbf{z}^{T}\textbf{H}\textbf{z}+ \textbf{q}_{k}^{T}\textbf{z}+{r}_{k}\\
 \text{subject to}
~&~ \textbf{F}\textbf{z} \leq \textbf{g}.
     \end{aligned}
\end{equation}
Note that here the parameters $\textbf{q}_{k},{r}_{k}$ and $\textbf{g}$ are functions of $\textbf{x}_{k}$. Therefore the current state information is required for solving this optimization problem.
\par Another way to reduce the online computation is to use a control horizon $N_C$ lesser than the prediction horizon $N.$ This in turn reduces the number of optimization variables. In this case, we define the control sequence as
$\textbf{U}_{k} = 
\big(\textbf{u}_{k|k},...,\textbf{u}_{k+{N}_C-1|k},\textbf{0},...,\textbf{0}\big)$ and this reduces the number of decision variables in $\textbf{z}$ to $mN_{c}$.
\newpage
\subsection{LMPC: Set point tracking}
So far we considered the stabilization problem in MPC for which the reference $\textbf{x}_{r}=0.$ In this section we discuss the set point tracking problem for which the reference $\textbf{x}_{r}\neq 0,$ and the objective is to track the nonzero set point. For the nonzero reference $\textbf{x}_{r},$ the steady state value of the control input will be nonzero, i.e. $\textbf{u}_{r}\neq 0$ and in steady state we have $\textbf{x}_{k+1}=\textbf{x}_{k}=\textbf{x}_{r}.$ Substituting this in (1) gives 
\begin{equation}
    \textbf{x}_{r}=\textbf{A}\textbf{x}_{r}+\textbf{B}\textbf{u}_{r}
    \implies  \textbf{u}_{r}=\textbf{B}^{-1}(\textbf{I}-\textbf{A})\textbf{x}_{r}
\end{equation}
where $\textbf{B}^{-1}$ is the pseudo-inverse.
The set point tracking can be transferred to a stabilization problem by defining the error state and control $\textbf{x}_{e_k}=\textbf{x}_{k}-\textbf{x}_{r},\textbf{u}_{e_k}=\textbf{u}_{k}-\textbf{u}_{r}$ and consider the error dynamics for MPC design which gives
\begin{equation}
    \begin{aligned}
    \textbf{x}_{e_{k+1}}&=\textbf{x}_{k+1}-\textbf{x}_{r}=\textbf{A}\textbf{x}_{k}+\textbf{B}\textbf{u}_{k}-\textbf{x}_{r}=\textbf{A}\textbf{x}_{k}-\textbf{A}\textbf{x}_{r}+\textbf{B}\textbf{u}_{k}-\textbf{x}_{r}+\textbf{A}\textbf{x}_{r}\\
    &=\textbf{A}[\textbf{x}_{k}-\textbf{x}_{r}]+\textbf{B}[\textbf{u}_{k}-\textbf{B}^{-1}(\textbf{I}-\textbf{A})\textbf{x}_{r}]=\textbf{A}\textbf{x}_{e_k}+\textbf{B}\textbf{u}_{e_k}.
    \end{aligned}
\end{equation}
Using the error state and control vectors the constraints can be rewritten as
\begin{equation}
\begin{aligned}
  & \textbf{F}_{\textbf{x}}\textbf{x} \leq \textbf{g}_{\textbf{x}}\implies \textbf{F}_{\textbf{x}}(\textbf{x}_{e_k}+\textbf{x}_{r}) \leq \textbf{g}_{\textbf{x}} \implies \textbf{F}_{\textbf{x}}\textbf{x}_{e_k} \leq \textbf{g}_{\textbf{x}}-\textbf{F}_{\textbf{x}}\textbf{x}_{r}  \\
  & \textbf{F}_{\textbf{u}}\textbf{u} \leq \textbf{g}_{\textbf{u}} \implies \textbf{F}_{\textbf{u}}(\textbf{u}_{e_k}+\textbf{u}_{r}) \leq \textbf{g}_{\textbf{u}} \implies  \textbf{F}_{\textbf{u}}\textbf{u}_{e_k} \leq \textbf{g}_{\textbf{u}}-\textbf{F}_{\textbf{u}}\textbf{u}_{r}.
  \end{aligned} 
\end{equation}
Now, the matrices $\textbf{F}_{\textbf{X}},\textbf{g}_{\textbf{X}},\textbf{F}_{\textbf{U}},\textbf{g}_{\textbf{U}}$ can be defined as in (12) in which $\textbf{g}_{\textbf{x}},\textbf{g}_{\textbf{u}}$ are replaced by $\textbf{g}_{x}-\textbf{F}_{\textbf{x}}\textbf{x}_{r},\textbf{g}_{u}-\textbf{F}_{\textbf{u}}\textbf{u}_{r}.$ We define $\textbf{z}=\left[\begin{matrix} \textbf{X}_{e_k}\\ \textbf{U}_{e_k} \end{matrix}\right]=\left[ \begin{matrix} \textbf{X}_{k}-\textbf{X}_{r} \\\textbf{U}_{k}-\textbf{U}_{r} \end{matrix}\right]$ and the optimization problem is obtained as in (15), solving which the optimal control input for the MPC problem is obtained as
\begin{equation}
    \textbf{u}_{k}=[\textbf{U}_{e_k}^{*}]_{1}+\textbf{u}_{r}.
\end{equation}

\subsection{LMPC: Numerical examples}
We consider an LTI system with system and input matrices as follows
\begin{equation}
    \textbf{A}=\left[\begin{matrix} 0.9 &0.2 \\-0.4 & 0.8
\end{matrix}\right] \hspace{1cm} \textbf{B}=\left[\begin{matrix} 0.1 \\0.01
\end{matrix}\right].
\end{equation}
 The simulation parameters are chosen as $N_{T}=50,$ $N=5,$ $\textbf{Q}=\textbf{I}_{2},\textbf{R}=1$ and $\textbf{x}_{0}=\left[\begin{matrix} 10 & 5
\end{matrix}\right]^{T}.$ The constraint set is defined as in (2) with 
\begin{equation}
    \textbf{F}_{\textbf{x}}=\left[\begin{matrix} 1 &0 \\0 & 1\\-1 &0 \\0&-1
\end{matrix}\right] \hspace{0.3cm} \textbf{g}_{\textbf{x}}=\left[\begin{matrix} 10 \\10\\10\\10
\end{matrix}\right]\hspace{0.3cm} \textbf{F}_{\textbf{u}}=\left[\begin{matrix} 1\\-1 \end{matrix} \right] \hspace{0.3cm} \textbf{g}_{\textbf{u}}=\left[\begin{matrix} 1\\1 \end{matrix} \right]
\end{equation}
which is equivalent to $-10 \leq x_{1k} \leq 10,$ $-10 \leq x_{2k} \leq 10,$ $-1 \leq u_{k} \leq 1.$
    The response of the LTI system with the MPC scheme is given in Fig. 2(a). The response shows the states converge to the origin and the constraints are satisfied. Similarly, for the set-point tracking problem, the state reference is chosen as $\textbf{x}_{r}=\left[\begin{matrix}3&2 \end{matrix} \right]^{T}$ for which the steady-state control input is obtained by solving (21) for the linear system (25) which gives $u_{r}=0.59$ which satisfies the control constraints. The simulation response for the set-point tracking is given in Fig. 2(b), which shows the state converges to the desired reference.
\begin{figure}[H]
 		\begin{center}
 		\includegraphics [scale=.5] {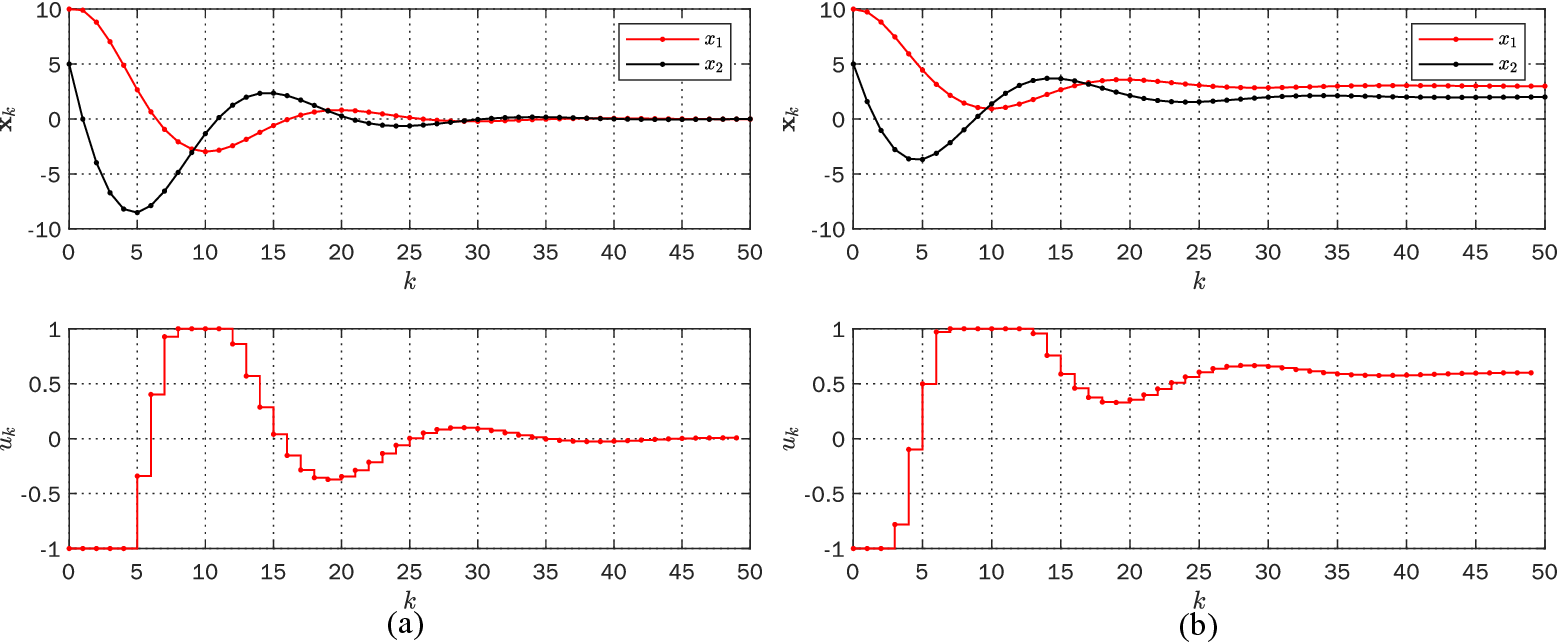}
 		\caption{{\footnotesize  LMPC response \hspace{.15cm} (a) Stabilization \hspace{.1cm}(b) Set point tracking.}}
 	\end{center}
 \end{figure}
\newpage
\section{MPC of Nonlinear Systems}
In this section, we discuss the basic concept and numerical implementation of Nonlinear MPC (NMPC).
\subsection{NMPC: Problem formulation}
Consider the discrete-time nonlinear system defined by the state equation:
\begin{equation}
    \textbf{x}_{k+1}=\textbf{f}(\textbf{x}_{k},\textbf{u}_{k})
\end{equation}
where $k\in \{0,1,...,N_{T}-1\}$ is the discrete time instant, $\textbf{x}_{k}\in \mathbb{X} \subseteq \mathbb{R}^{n}$ is the state vector,  $\textbf{u}_{k}\in \mathbb{U} \subseteq \mathbb{R}^{m}$ is the input vector and $\textbf{f}:\mathbb{X} \times \mathbb{U} \rightarrow \mathbb{X}$ is the nonlinear mapping which maps the current state $\textbf{x}_{k}$ to the next state $\textbf{x}_{k+1}$ under the control action $\textbf{u}_{k}.$
The constraint sets $\mathbb{X}$ and $\mathbb{U}$ are defined as in (2) and the cost function is chosen as a quadratic function as in (5). Then the MPC problem for nonlinear systems is defined as follows:
\begin{problem}
For the nonlinear system (27) with the current state $\textbf{x}_{k|k}=\textbf{x}_{k}$,
compute the control sequence $\textbf{U}_{k}$ by solving the optimization problem
\begin{equation}
    \begin{aligned}
    \underset{\textbf{U}_{k}}{\mbox{\text{inf}}} \hspace{.2cm} &J_{k}\\ \text{subject to}
    ~&~ \textbf{U}_{k} \in \mathbb{U}^{N}, \hspace{.2cm} \textbf{X}_{k}\in \mathbb{X}^{N+1},\hspace{.4cm}k\in \mathbb{T}\\
    ~&~ \textbf{x}_{i+1|k}=\textbf{f}(\textbf{x}_{i|k},\textbf{u}_{i|k}), \hspace{.3cm}k\in \mathbb{T},i=k,...,k+N-1.
    \end{aligned}
\end{equation}
\end{problem}.
\subsection{NMPC: Algorithm}
\par By defining  $\textbf{X}_{k}$ and 
$\textbf{U}_{k}$ as in (8) we can rewrite the cost function and constraints for the nonlinear MPC problem as 
\begin{equation}
    J_{k}=  \textbf{X}_{k}^{T}\textbf{Q}_{\textbf{X}} \textbf{X}_{k}+\textbf{U}_{k}^{T}\textbf{R}_{\textbf{U}} \textbf{U}_{k}
\end{equation}
and
\begin{equation}
\begin{aligned}
  & \textbf{F}_{\textbf{X}}\textbf{X}_{k} \leq \textbf{g}_{\textbf{X}}\\
   & \textbf{F}_{\textbf{U}}\textbf{U}_{k} \leq \textbf{g}_{\textbf{U}}\\ 
&\textbf{f}_{eq}(\textbf{X}_{k},\textbf{U}_{k})=0
    \end{aligned}
\end{equation}
where 
\begin{equation}
    \textbf{f}_{eq}(\textbf{X}_{k},\textbf{U}_{k})= \left[\begin{matrix}  \textbf{x}_{k|k}-\textbf{x}_{k}\\\textbf{x}_{k+1|k}-\textbf{f}(\textbf{x}_{k|k},\textbf{u}_{k|k})\\ \vdots \\
    \textbf{x}_{k+N|k}-\textbf{f}(\textbf{x}_{k+{N}-1|k},\textbf{u}_{k+{N}-1|k})
\end{matrix}\right].
\end{equation}
Now, by defining $\textbf{z},\textbf{H},\textbf{F},\textbf{g}$ as in (14) the optimization problem is represented as a nonlinear programming problem as below
\begin{equation}
\begin{aligned}
 \underset{\textbf{z}}{\text{inf}} ~&~     \textbf{z}^{T}\textbf{H}\textbf{z}\\
 \text{subject to}
~&~ \textbf{F}\textbf{z} \leq \textbf{g}\\
  ~&~ \textbf{f}_{eq}(\textbf{z}) = 0.
     \end{aligned}
\end{equation}
Here the equality constraint is nonlinear which makes the optimization problem a nonlinear programming problem.
In MPC this optimization problem is solved during every time instant $k$ and the first element of $\textbf{U}_{k}^{*}$ is applied to the system, i.e., the control input with MPC  is 
\begin{equation}
    \textbf{u}_{k}=[\textbf{U}_{k}^{*}]_{1}=\textbf{u}_{k|k}^{*}.
\end{equation}
The algorithm for nonlinear MPC is summarized below:
\begin{algorithm}[H]
 \small
	\begin{algorithmic}[1] 
	
	\STATE Require $\textbf{f},N_{T},N,n,m,\textbf{Q},\textbf{R},\textbf{Q}_{N_T},\textbf{F}_{\textbf{x}},\textbf{g}_{\textbf{x}},\textbf{F}_{\textbf{u}},\textbf{g}_{\textbf{u}}$
		\STATE Initialize $\textbf{x}_{0},\textbf{z}_{0}$
		\STATE Construct $\textbf{Q}_{X},\textbf{R}_{\textbf{U}},\textbf{H},\textbf{F},\textbf{g}$
		\FOR  {$k= 0~to~ N_{T}-1 $}
		\STATE $\textbf{x}_{k}=[\textbf{X}]_{k+1}$ (obtain $\textbf{x}_{k}$ from measurement or estimation)
		\STATE  Compute $\textbf{z}^{*}=\left[\begin{matrix}\textbf{X}_{k}^{*}\\\textbf{U}_{k}^{*}\end{matrix}\right]$ by solving the optimization problem (32)
		\STATE Apply  $\textbf{u}_{k}=[\textbf{U}_{k}^{*}]_{1}$ to the system
		\STATE Update $\textbf{z}_{0}=\textbf{z}^{*}$
		\ENDFOR
	\end{algorithmic}
	\caption{: NMPC}
\end{algorithm}
The optimization problem (32) can be solved using the MATLAB function for solving constrained optimization problems:
    \begin{equation}
        \textbf{z}^{*}=\mbox{fmincon}(f,\textbf{z}_{0},\textbf{F},\textbf{g},\textbf{lb},\textbf{ub},\textbf{f}_{eq}).
    \end{equation}
\subsection{NMPC: Set point tracking}
Here we discuss the set point tracking problem for nonlinear systems for which the reference $\textbf{x}_{r}\neq 0$. The reference value or steady state value of the control input $\textbf{u}_{r}$ is computed by solving the steady-state equation
\begin{equation}
    \textbf{x}_{r}=\textbf{f}(\textbf{x}_{r},\textbf{u}_{r}).
\end{equation}
By defining the error state and control vector as $\textbf{x}_{e_k}=\textbf{x}_{k}-\textbf{x}_{r}, \textbf{u}_{e_k}=\textbf{u}_{k}-\textbf{u}_{r},$ the constraints can be rewritten as in (23). Similarly, the equality constraint  becomes
\begin{equation}
    \textbf{x}_{e_{k+1}}=\textbf{f}(\textbf{x}_{k},\textbf{u}_{k})-\textbf{x}_{r}=\textbf{f}(\textbf{x}_{e_k}+\textbf{x}_{r},\textbf{u}_{e_k}+\textbf{u}_{r})-\textbf{x}_{r}.
\end{equation}
Now by defining $\textbf{z}=\left[\begin{matrix} \textbf{X}_{e_k}\\ \textbf{U}_{e_k} \end{matrix}\right]=\left[ \begin{matrix} \textbf{X}_{k}-\textbf{X}_{r} \\\textbf{U}_{k}-\textbf{U}_{r} \end{matrix}\right],$ the optimization problem is obtained as in (32), solving which the optimal control input for the MPC problem is obtained as
\begin{equation}
    \textbf{u}_{k}=[\textbf{U}_{e_k}^{*}]_{1}+\textbf{u}_{r}.
\end{equation}
\subsection{NMPC: Numerical examples}
We consider the discrete-time model of the simple pendulum system which is defined by the state equation:
\begin{equation}
\textbf{x}_{k+1}=\textbf{f}(\textbf{x}_{k},\textbf{u}_{k})=\left[\begin{matrix}x_{1_k}+T x_{2_k}\\x_{2_k}+T\big(-\frac{g}{l}sin(x_{1_k})-\frac{B}{Ml^{2}}x_{2_k}+\frac{1}{Ml^{2}}u_k\big) \end{matrix}\right]
\end{equation}
where $M$ is the mass of the simple pendulum, $B$ is the friction coefficient, $l$ is the length of the pendulum, $g$ is the acceleration due to gravity and $T$ is the sampling time. The system parameters are chosen as $M=1,B=1,l=1,g=9.8,T=0.1$ and simulation parameters are chosen as $\textbf{Q}=\textbf{I}_{2},\textbf{R}=1$ and $\textbf{x}_{0}=\left[\begin{matrix} 2 & 1
\end{matrix}\right]^{T}.$ The constraint set parameters is defined as 
\begin{equation}
    \textbf{F}_{\textbf{x}}=\left[\begin{matrix} 1 &0 \\0 & 1\\-1 &0 \\0&-1
\end{matrix}\right] \hspace{0.3cm} \textbf{g}_{\textbf{x}}=\left[\begin{matrix} 5 \\5\\5\\5
\end{matrix}\right]\hspace{0.3cm} \textbf{F}_{\textbf{u}}=\left[\begin{matrix} 1\\-1 \end{matrix} \right] \hspace{0.3cm} \textbf{g}_{\textbf{u}}=\left[\begin{matrix} 0.1\\0 \end{matrix} \right]
\end{equation}
which is equivalent to $-5 \leq x_{1k} \leq 5,$ $-5 \leq x_{2k} \leq 5,$ $0 \leq u_{k} \leq 0.1.$
    The response of the simple pendulum with the MPC scheme is given in Fig. 3(a). The response shows the states converge to the origin and the constraints are satisfied. Similarly for the set-point tracking problem, the state reference is chosen as $\textbf{x}_{r}=\left[\begin{matrix}0.5&0 \end{matrix} \right]^{T}$ for which the steady-state control input is obtained by solving (36) for the nonlinear system (38) which gives $u_{r}=Mgl\hspace{.1cm} sin(x_{1_r})=4.69.$ Hence we set the maximum value of control input as $5$ for the set-point tracking problem. The simulation response for the set-point tracking is given in Fig. 3(b) which shows the state converges to the desired reference.
\begin{figure}[H] 
 		\begin{center}
 		\includegraphics [scale=.5] {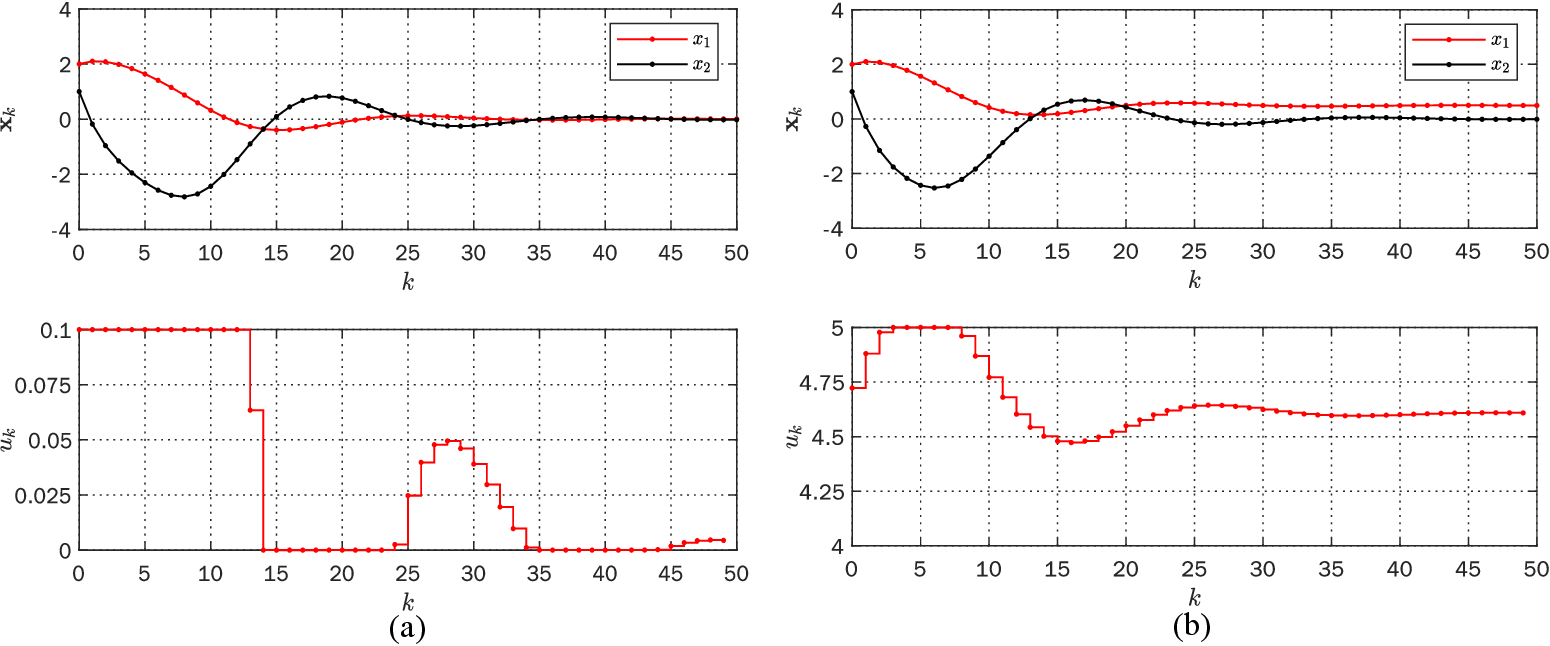}
 		\caption{{\footnotesize  NMPC response \hspace{.15cm} (a) Stabilization \hspace{.1cm}(b) Set point tracking.}}	
 	\end{center}
 \end{figure}
\newpage
\section{Feasibility,  Stability, and Optimality}
In this section, we study the feasibility, stability, and optimality of the MPC scheme. We start with the feasibility which deals with the existence of the optimal solution. The MPC problem is feasible, if there exists an optimal solution $\textbf{z}^{*}$ for the optimization problem at each time instant that satisfies all the constraints.
Whenever there are constraints on states the optimization problem becomes more complicated. In that case, we have to select the control sequence $\textbf{U}_{k}$ in such a way that the corresponding predicted state sequence $\textbf{X}_{k}$ does not violate the state constraints. This leads to the idea of feasibility and feasible sets. We denote $\mathbb{U}_{fk} \subseteq \mathbb{U}^{N}$ as the feasible set of control inputs
\begin{equation}
 \mathbb{U}_{fk}=\{   \textbf{U}_{k}\in \mathbb{U}^{N}: \textbf{X}_{k}(\textbf{x}_{k},\textbf{U}_{k}) \in \mathbb{X}^{{N}+1} \}.
\end{equation}  
Note that $\mathbb{U}_{fk}$ depends on the current state $\textbf{x}_{k},$ i.e., $\mathbb{U}_{fk}=\mathbb{U}_{f}(\textbf{x}_{k})$ and we denote it as $\mathbb{U}_{fk}$ to simplify notations. It also depends on the prediction horizon $N$ which we are considering as fixed here.
The number of elements in $\mathbb{U}_{fk}$ decreases when $\textbf{x}_{k}$ is closer to the boundary of $\mathbb{X}$ and when $\textbf{x}_{k}$
moves away from the boundary more and more $\textbf{U}_{k}$ becomes feasible and when $\textbf{x}_{k}$ is sufficiently far from the boundary we have $\mathbb{U}_{fk}=\mathbb{U}^{N},$ i.e., all the control sequences are feasible,  which is same as the unconstrained state case. This situation is demonstrated in Fig. 4 in which for Fig. 4(a) the current state $\textbf{x}_{k}$ is closer to the boundary of the constraint set. In this case, the predicted state sequences 2 and 3 violate the state constraints, hence the corresponding control sequences will not be feasible. But, for Fig. 4(b) the current state $\textbf{x}_{k}$ is sufficiently far away from the boundary of the constraint set which makes all the 3 predicted state sequences to stay within the constraint set. Consequently, all the 3 control sequences will be feasible. \par The MPC problem is said to be feasible for $\textbf{x}_{k}\in \mathbb{X}$ if $\mathbb{U}_{fk}$ is nonempty. This also ensures the existence of a solution to the optimization problem. Clearly, for the unconstrained state case, the MPC problem is always feasible, and for the constrained state case, feasibility depends on the current state. We denote the set for feasible states by $\mathbb{X}_{fk} \subseteq \mathbb{X}$ which is defined as
\begin{equation}
 \mathbb{X}_{fk}=\{   \textbf{x}_{k}\in \mathbb{X}: \mathbb{U}_{fk} \neq \phi \}.
\end{equation} 
In general, if $\mathbb{X}_{fk}$ and $\mathbb{U}_{fk}$ are the feasible set of states and control sequences during time instant $k.$ Then the MPC control law is computed by solving the optimization problem:
\begin{equation}
 \begin{aligned}
\underset{\textbf{U}_{k}\in \mathbb{U}_{f_k}}{\text{inf}}  ~&~ J_{k}(\textbf{x}_{k},\textbf{U}_{k})\\
\text{subject to}
  ~&~ \textbf{x}_{i+1|k}=\textbf{f}(\textbf{x}_{i|k},\textbf{u}_{i|k}), \hspace{.3cm}k\in \mathbb{T}, i=k,...,k+N-1.
     \end{aligned}
\end{equation}
Clearly, every control sequence in the set $\mathbb{U}_{f_k}$ results in a predicted state sequence that satisfies the state constraints. Hence there is no need to include the state constraints in the optimization problem here. 
The notation $\mathbb{X}_{fk}$ is more general and covers the time-varying systems also, and for time-invariant systems, the index $k$ can be omitted which simplifies the notation to $\mathbb{X}_{f}.$
\par Another important concept associated with feasibility is persistent feasibility.
The MPC problem is said to be persistently feasible, if the feasibility of initial state $\textbf{x}_{0}$ guarantees the feasibility of future states $\textbf{x}_{k}, k=1,2,...,N_T$ under the dynamics, i.e. $\mathbb{U}_{f0}\neq \phi \implies \mathbb{U}_{fk} \neq \phi,$ $\forall k=1,2,...,N_{T}$. Persistent feasibility depends on the system dynamics, prediction horizon $N,$ and the constrained sets $\mathbb{X},\mathbb{U}$.

\begin{figure}[H]
 		\begin{center}
 		\includegraphics [scale=.4] {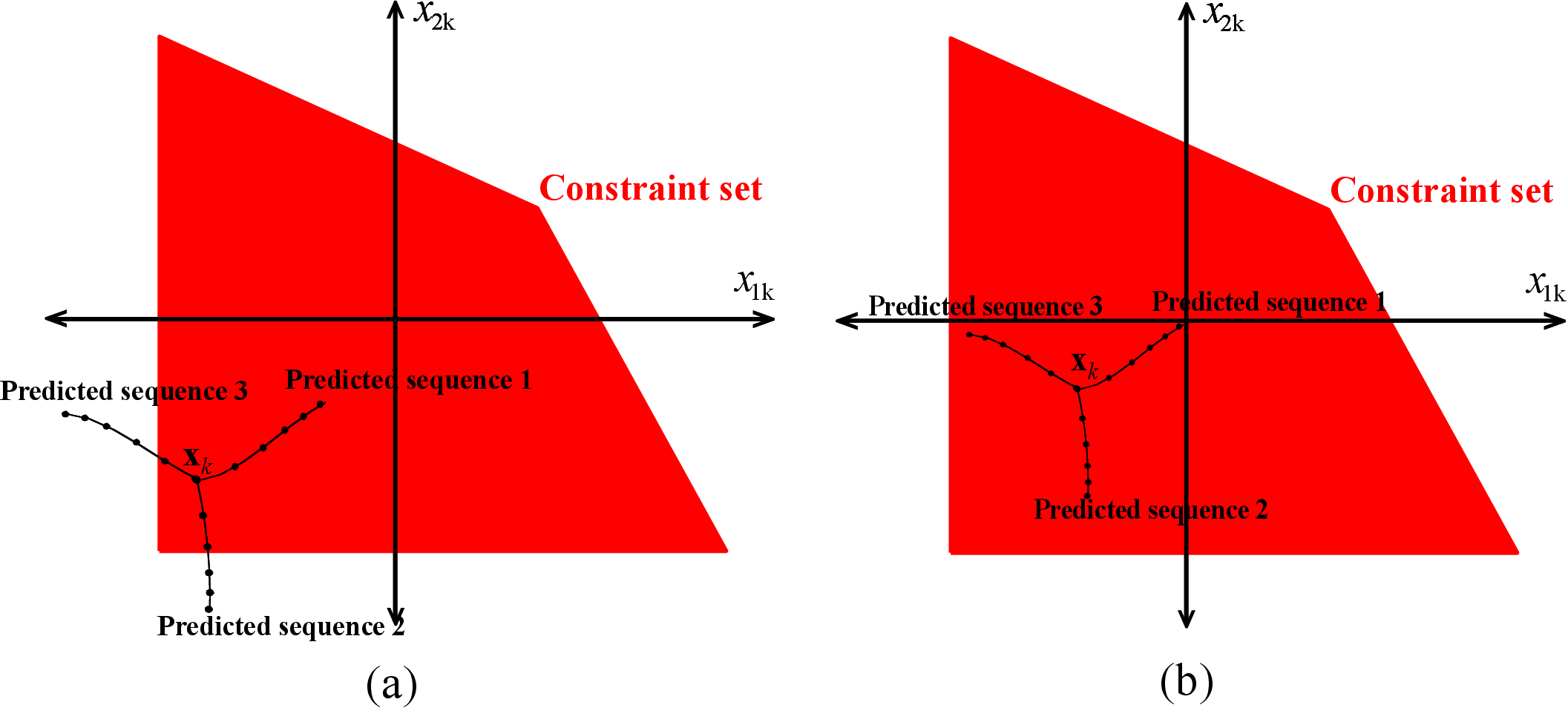}
 		\caption{{\footnotesize  Feasibility \hspace{.15cm} (a) State near the boundary of $\mathbb{X}$ \hspace{.1cm}(b) State away from the boundary of $\mathbb{X}$.}}	
 	\end{center}
 \end{figure}
\newpage
\par Next, we discuss stability, which deals with the convergence of the solution, i.e., whether the state trajectory under the MPC scheme converges to the desired reference or equilibrium point. In MPC, the stability analysis is mainly based on the Lyapunov approach in which the basic idea is to design the control scheme in such a way that the optimal cost function becomes a Lyapunov function, i.e., $V_{k}=J_{k}^{*}$ and it satisfies
\begin{equation}
    \varDelta V=J_{k+1}^{*}(\textbf{x}_{k+1})-J_{k}^{*}(\textbf{x}_{k})<0.
\end{equation}
There exist different variants of the criteria (43) which give a different upper bound for $\varDelta V.$ In general, for stabilizable LTI systems, by properly selecting the terminal weighting matrix and constraints, the value function of the MPC scheme can be made as a Lyapunov function. However, this may not always be possible for nonlinear systems. The terminal weighting matrix $\textbf{Q}_{N}$ and terminal constraints $\textbf{F}_{\textbf{x}_{N}},\textbf{g}_{\textbf{x}_{N}}$ can be easily incorporated in the MPC algorithm by adding them in $\textbf{Q}_{\textbf{X}},\textbf{F}_{\textbf{X}},\textbf{g}_{\textbf{X}}$ which results in:
\begin{equation}
\textbf{Q}_{\textbf{X}} = \left[\begin{matrix}
\textbf{Q} & \dots & \textbf{0} & \textbf{0}\\ \vdots &  &\vdots & \vdots \\ \textbf{0} & \dots & \textbf{Q} & \textbf{0} \\ \textbf{0}& \dots  &\textbf{0} & \textbf{Q}_{N}
\end{matrix}\right], \hspace{.2cm}\textbf{F}_{\textbf{X}} = \left[\begin{matrix}
\textbf{F}_{\textbf{x}} & \dots & \textbf{0} & \textbf{0}\\ \vdots &  &\vdots & \vdots \\ \textbf{0} & \dots & \textbf{F}_{\textbf{x}} & \textbf{0} \\ \textbf{0}& \dots  &\textbf{0} & \textbf{F}_{\textbf{x}_{N}}
\end{matrix}\right], \hspace{.2cm}\textbf{g}_{\textbf{X}} = \left[\begin{matrix}
\textbf{g}_{\textbf{x}}\\\vdots\\\textbf{g}_{\textbf{x}} \\\textbf{g}_{\textbf{x}_{N}}
\end{matrix}\right].
\end{equation}
\par Finally, optimality is a term associated with the performance of the solution and depends on how fast the trajectory converges to the equilibrium point and how much control effort is required.
When it comes to optimality, the MPC schemes usually result in a suboptimal solution. This is because of the reason that, in MPC during every time instant we optimize the performance over the prediction horizon, not the entire time horizon. Consequently, as the prediction horizon increases the MPC control law becomes more optimal and, in general as $N \rightarrow N_T$, the control law becomes optimal.
\section{Further Reading}
This tutorial attempts to discuss the basic theory of MPC and its numerical implementation using MATLAB. For a more detailed study on linear MPC and nonlinear MPC, one can refer to \cite{bFA} and \cite{bLJ}, respectively. For a better understanding of the numerical optimization methods for solving linear programming, quadratic programming, and nonlinear programming problems, one may refer to \cite{bDY}.
For related and advanced topics in MPC such as LQR, Kalman filter,   adaptive MPC, Robust MPC, and Distributed MPC, one can refer \cite{bDM}-\cite{bMA1}.
A lecture series based on this tutorial can be found at \cite{bMA2}.
The MATLAB codes for the MPC examples discussed in this paper are available at \cite{bMA3}.
\bibliographystyle{unsrt}

\end{document}